\documentclass[11pt]{article}          
\makeatletter
\def\@begintheorem#1#2{\trivlist%
 \item[\hskip \labelsep{\sffamily\bfseries #2\ #1}]\itshape}
\newtheorem{teo}{Theorem}[section]
\newtheorem{defi}[teo]{Definition}
\newtheorem{cor}[teo]{Corollary}
\newtheorem{lem}[teo]{Lemma}
\newtheorem{pro}[teo]{Proposition}
\newtheorem{_rem}[teo]{Remark}
\newtheorem{_eje}[teo]{Example}
\newenvironment{rem}{\def\@begintheorem##1##2{\trivlist%
 \item[\hskip\labelsep{\sffamily\bfseries ##2\ ##1}]}\begin{_rem}}{\end{_rem}}

\makeatother

\newenvironment{beweis}{{\em Proof:}}{\hfill $\rule{2mm}{2mm}$
\vspace{3mm}

}

\DeclareMathAlphabet{\Ma}{U}{msa}{m}{n}
\DeclareMathAlphabet{\Mb}{U}{msb}{m}{n}
\DeclareMathAlphabet{\Meuf}{U}{euf}{m}{n}

\def\got#1{\Meuf{#1}}

\DeclareSymbolFont{ASMa}{U}{msa}{m}{n}
\DeclareSymbolFont{ASMb}{U}{msb}{m}{n}
\DeclareMathSymbol{\hrist}{\mathord}{ASMa}{"16}
\DeclareMathSymbol{\varkappa}{\mathalpha}{ASMb}{"7B}
\DeclareMathSymbol{\CrPr}{\mathord}{ASMb}{"6F}

\newfont{\EinsFont}{cmr7 scaled 1070}
\def\EINS{{\mathchoice{
 \mbox{\unitlength1cm\begin{picture}(.25,.2)\put(0,0){$1$}%
 \put(0.105,0){{\mbox{\fontfamily{cmr}\upshape\small l}}}\end{picture}}}{%
 \mbox{\unitlength1cm\begin{picture}(.25,.2)\put(0,0){$1$}%
 \put(0.105,0){{\mbox{\fontfamily{cmr}\upshape\small l}}}\end{picture}}}{%
 \mbox{\unitlength1cm\begin{picture}(.18,.15)\put(0,0){$\scriptstyle 1$}%
 \put(0.07,0){{\mbox{\fontfamily{cmr}\upshape\EinsFont l}}}\end{picture}}}{%
 \mbox{\unitlength1cm\begin{picture}(.18,.15)\put(0,0){$\scriptstyle 1$}%
 \put(0.07,0){{\mbox{\fontfamily{cmr}\upshape\EinsFont l}}}\end{picture}}}}}

\def\restriction{{\mathchoice{
 \mbox{\unitlength1cm\begin{picture}(.2,.4)%
  \bezier{5}(.07,.3)(.1,.27)(.13,.24)%
  \put(.07,.35){\line(0,-1){.5}}\end{picture}}}{
 \mbox{\unitlength1cm\begin{picture}(.2,.4)%
  \bezier{5}(.07,.3)(.1,.27)(.13,.24)%
  \put(.07,.35){\line(0,-1){.5}}\end{picture}}}{
  \hrist}{\hrist}}}

  \def\al #1.{{\mathcal{#1}}}
  \def\ot #1.{{\got{#1}}}
  \def\C{\Mb{C}}
  
  \def\Z{\Mb{Z}}
  \def\T{\Mb{T}}

\title{\bf Dual group actions on C*--algebras and 
        their description by Hilbert extensions}
\author{
 {\sc Hellmut Baumg\"artel}  \\[2mm] 
 {\footnotesize Mathematical Institute, University of Potsdam,}     \\   
 {\footnotesize Am Neuen Palais 10, Postfach 601~553,}             \\ 
 {\footnotesize D--14415 Potsdam, Germany.}                        \\
 {\footnotesize baumg@rz.uni-potsdam.de}\\
 {\footnotesize FAX: 49--(0)331--9771299}
\and    
 {\sc Fernando Lled\'o }\thanks{On leave from {\em
  Mathematical Institute, University of Potsdam,
  Am Neuen Palais 10, Postfach 601~553,
  D--14415 Potsdam, Germany.} }  \\[2mm] 
 {\footnotesize Max--Planck--Institut f\"ur Gravitationsphysik,}     \\   
 {\footnotesize Albert--Einstein--Institut,}
 {\footnotesize Am M\"uhlenberg 1,}             \\ 
 {\footnotesize D--14476 Golm, Germany.}                        \\
 {\footnotesize lledo@aei-potsdam.mpg.de}}

\date{\today{}}

\begin{document}
\maketitle
\begin{abstract}
Given a $C^{\ast}$-algebra ${\cal A}$, a discrete abelian group
${\cal X}$ and a homomorphism $\Theta\colon\ {\cal X}\rightarrow
\mbox{Out}\,{\cal A},$ defining the dual action group
$\Gamma\subset \mbox{aut}\,{\cal A}$, the paper contains 
results on existence and characterization of Hilbert extensions of
$\{\al A.,\Gamma\}$, where the action is given by $\hat{\cal X}$.
They are stated at the (abstract) C*--level and can therefore 
be considered as a refinement of the extension results given
for von Neumann algebras for example by Jones \cite{Jones80} or 
Sutherland \cite{Sutherland80,SutherlandIn82}. A Hilbert extension
exists iff there is a generalized 2--cocycle. These results generalize
those in \cite{DHR69b}, which are formulated 
in the context of superselection theory, where it is assumed that the 
algebra $\al A.$ has a trivial center, i.e.~$\al Z.=\C\EINS$.
In particular the well--known ``outer characterization'' 
of the second cohomology 
$H^{2}({\cal X},{\cal U}({\cal Z}),\alpha_{\cal X})$ can be reformulated:
there is a bijection to the set of all
${\cal A}$--module isomorphy classes of Hilbert extensions.
Finally, a Hilbert space representation (due to Sutherland
\cite{Sutherland80,SutherlandIn82} in the von Neumann case)
is mentioned. The C*--norm of the Hilbert extension is expressed
in terms of the norm of this representation and it is linked to
the so--called regular representation appearing in superselection
theory.
\end{abstract}


\section{Introduction}
In the Doplicher/Roberts theory (e.g.~\cite{Doplicher89a,Doplicher90}) 
it is a central assumption that the center of the C*--algebra $\al A.$
with which one starts the analysis is trivial, 
i.e.~$\al Z.=\al Z.(\al A.)= \C\EINS$ (in a more categorial notation
the assumption reads $(\iota,\iota)=\C\EINS$, where $\iota$ denotes the
unit object of the strict monoidal C*--category, cf.~\cite{Doplicher89b}).
From a systematical point of view it is interesting to study the
properties and structural modifications of this theory if one assumes
the presence of a nontrivial center $\al Z.\supset\C\EINS$. For example,
if $(\al F.,\alpha_\al G.)$ is a Hilbert C*--system for a compact 
group $\al G.$ and if the corresponding fixed point algebra $\al A.$ 
has a nontrivial center that satisfies the relation $\al A.'\cap\al F.
=\al Z.$, then the Galois correspondence does not hold anymore, i.e.~we
have the proper inclusion $\alpha_\al G.\subset\mathrm{stab}\,\al A.$
in $\mathrm{aut}\,\al F.$ (cf.~\cite[Section~7]{Lledo97b}). Recall, that
in the trivial center situation it is a fundamental result of the
theory that $\alpha_\al G.=\mathrm{stab}\,\al A.$. As a further 
justification we can also mention that in other generalizations of the
Doplicher/Roberts theory as well as in some applications in mathematical
physics a nontrivial center plays, to a certain extent, a distinguished
role \cite{Longo97,pVasselli00,Fredenhagen92}. 

In the present paper we continue the analysis of the presence of a
nontrivial center in the construction of an extension algebra 
$\al F.$ (cf.~\cite{BaumgaertelIn99,BaumgaertelIn00}). 
In particular, we study what we call dual group actions
in the simple case where the group $\al X.$ is discrete and
abelian (cf.~with \cite{DHR69b} in the special case where
$\al Z.=\C\EINS$). This investigations will be done at the
abstract C*--level which is the context of the 
Doplicher/Roberts theory mentioned above (cf.~also \cite{Baumgaertel97}).
On the other hand the results can be considered as a refinement
of the study of twisted group algebras (twisted crossed products)
on the concrete von Neumann algebra level
(see e.g.~\cite{Busby70,Jones80,Sutherland80,SutherlandIn82}).
For example, the decisive C*--norm for the extension is defined 
intrisically and the natural representation (discussed e.g.~by
Sutherland) is related to the so--called regular representation
that appears in the superselection theory \cite{bBaumgaertel95}.
We hope that the present analysis will be useful to obtain a more general
`inversion' theorem, where endomorphisms of $\al A.$ are involved.
Indeed, the main theorems in Section~\ref{Sec3} suggest that for 
a more general inversion theory in the nontrivial center situation
the cohomological aspects may be essential.

The paper is structured in 5 sections: in the following section we will
introduce the notion of a Hilbert C*--system and study some 
properties of the group homomorphism $\Theta\colon\al X.\to\mathrm{Out}\,\al
A.$. Hilbert C*--systems are the result of the extension procedure
mentioned above. In Section~\ref{Sec3} we begin the study of the 
inverse (extension) problem: in particular it contains the result that
a Hilbert extension exists iff there is a generalized 
2--cocycle (to be defined there), and that in this case 
the set of all Hilbert extensions can be described in terms of the set
of center--valued 2--cocycles of 
$H^{2}({\cal X},{\cal U}({\cal Z}),\alpha_{\cal X})$
(cf.~Theorems~\ref{Main1} and \ref{Teo2}). In the next section
we relate the previously obtained results to the special case of the
Doplicher/Roberts frame, where $\al Z.=\C\EINS$. Finally, in
Section~\ref{Sec5} we give a representation of the Hilbert 
extension, which was already studied by Sutherland
\cite{Sutherland80,SutherlandIn82} in the von Neumann case. In particular,
we show that 
if there is a faithful state of $\al A.$, this representation
coincides with the so--called regular representation
that appears in superselection theory (cf.~e.g.~\cite{bBaumgaertel95})
and the intrinsic C*--norm turns out to be the operator norm of this
representation.

\section{Hilbert C*--systems}\label{Sec2}

A C*-algebra ${\cal F}$ together with a pointwise
norm-continuous group homomorphism ${\cal G}\ni g\rightarrow
\alpha_{g}\in \mbox{aut}\,{\cal F}$ of a locally compact group
${\cal G}$ is called a C*-system $\{{\cal F},\alpha_{\cal G}\}.$
Let ${\cal A}\subseteq{\cal F}$ be its fixed point algebra, i.e.~${\cal
A}:= \{A\in {\cal F}\mid\alpha_{g}A=A, g\in {\cal G}\}$. We denote
by ${\cal A}^{c}:={\cal F}\cap {\cal A}' \subseteq {\cal F}$ the
relative commutant of ${\cal A}$ w.r.t. ${\cal F}$. As is
well-known, $\alpha_{g}\restriction {\cal A}^{c}$ is an
automorphism of ${\cal A}^{c}$, so $\{{\cal A}^{c},\alpha_{\cal G}\}$ 
is also a C*-system. We call it the {\em assigned} C*-system.
The center ${\cal Z}({\cal A})$ is denoted by ${\cal Z}$.

In the following let ${\cal G}$ be compact and abelian so that
$\hat{\cal G}=:{\cal X}$ is abelian and discrete. The
corresponding spectral projections w.r.t. $\{{\cal
F},\alpha_{\cal G}\}$ are denoted by $\Pi_{\chi},\,\chi\in {\cal
X}$. Note that $\Pi_{\iota}{\cal F}={\cal A}$, where $\iota$
is the unit element of ${\cal X}$.

\begin{defi}
A C*-system $\{{\cal F},\alpha_{\cal G}\},\,{\cal
G}$ compact abelian, is called a {\em Hilbert} C*-{\em system}
if spec$\,\alpha_{\cal G}={\cal X}$ and if each spectral
subspace $\Pi_{\chi}{\cal F}$ contains a unitary $U_{\chi}$,
i.e. ${\cal U}(\Pi_{\chi}{\cal F})\neq \emptyset$.
\end{defi}

If $\{{\cal F},\alpha_{\cal G}\}$ is Hilbert, then
$\beta_{\chi}:= \mbox{ad}\,U_{\chi}\restriction {\cal A}$ is an
automorphism of ${\cal A}$, i.e.~$\beta_{\chi}\in \mbox{aut}\,{\cal A}.$
We denote by $\pi$ the canonical homomorphism of aut$\,{\cal A}$
onto Out$\,{\cal A}:=\mbox{aut}\,{\cal A}/\mbox{int}\,{\cal A}$,
where int$\,{\cal A}$ denotes the normal subgroup of all inner
automorphisms of ${\cal A}$. Then
\begin{equation}
{\cal X}\ni \chi\rightarrow \Theta(\chi):=\pi(\beta_{\chi})\in
\mbox{ Out}\,{\cal A}
\end{equation}
is a group homomorphism of ${\cal X}$ into Out$\,{\cal A}$, i.e.
we have

\begin{lem}\label{Lem1}  
To each Hilbert C*--system $\{{\cal F},\alpha_{\cal G}\}$, where
${\cal G}$ is compact abelian, there is canonically assigned a 
group homomorphism
$\Theta\colon\ {\cal X}\rightarrow \mbox{Out}\,{\cal A}$ given by (1).
\end{lem}
\begin{beweis}
Note that for $\chi_1,\chi_2\in\al X.$ we have that
$U_{\chi_1\chi_2}U_{\chi_2}^*
U_{\chi_1}^*\in\al A.$ and this implies that $\beta_{\chi_1\chi_2}
\circ\beta_{\chi_2}^{-1}\circ\beta_{\chi_1}^{-1}\in\mathrm{int}\,\al A.$.
\end{beweis}

We mention next the characterization of those Hilbert C*-systems where
$\Theta$ is an isomorphism and of those where the classes
$\Theta(\chi)$ are pairwise disjoint. Recall that $\alpha,\beta\in
\mbox{aut}\,{\cal A}$ are called disjoint if 
\[
  (\alpha,\beta):= \{X\in\al A.\mid 
                   X\alpha (A)=\beta(A)X \quad\mathrm{for~all}
                   \quad A\in {\cal A}\} =0\,.
\]

\begin{pro}
\begin{itemize}
\item[(i)] $\Theta$ is a monomorphism iff no spectral
subspace $\Pi_{\chi}{\cal A}^{c},\chi\neq\iota$, of the
assigned C*--system contains a unitary.
\item[(ii)] The classes $\Theta(\chi)$ are pairwise disjoint iff ${\cal
A}^{c} = {\cal Z}$, i.e. the relative commutant
coincides with the center of ${\cal A}$.
\end{itemize}
\end{pro}
\begin{beweis}
For one of the directions of part (i) take 
a unitary $U_\chi\in \Pi_\chi (\al
A.^c)$ with $\iota\not=\chi\in\al X.$, so that the corresponding
$\beta_\chi=\mathrm{id}$ and $\pi(\beta_\chi)=\mathrm{int}\,\al A.$. 
Thus $\Theta$ is not injective. For the other implication take
$\al X.\ni\chi_0\not=\iota$ with $\chi_0\in\mathrm{ker}\,\Theta$,
i.e.~$\Theta(\chi_0)=\mathrm{int}\,\al A.$. Thus there exists a unitary
$V\in\al U.(\al A.)$ with $\mathrm{ad}\,V=\mathrm{ad} \,U_{\chi_0}$.
From this we get 
$V^*U_{\chi_0}\in\al U.(\al A.^c)\cap\Pi_{\chi_0}(\al F.)$, 
i.e.~$\Pi_{\chi_0}{\cal A}^{c}\neq\emptyset$. 

Finally, part (ii) follows from \cite[Lemma~10.1.8]{bBaumgaertel92}.
\end{beweis}

We mention several useful concepts for Hilbert C*-systems
$\{{\cal F},\alpha_{\cal G}\}$ with a compact abelian group.

\begin{defi}
$\beta\in\mbox{aut}\,{\cal A}$ is called a {\em
canonical automorphism} if $\beta:= \mbox{ad}\,V\restriction
{\cal A}\,,V\in \bigcup_{\chi\in{\cal X}}{\cal
U}(\Pi_{\chi}{\cal F}).$ The set of all canonical automorphisms
is denoted by $\Gamma$.
\end{defi}
\begin{rem}
Note that for the set of canonical automorphisms we have
$\mathrm{int}\,\al A.\subseteq\Gamma\subseteq\mathrm{aut}\,\al A.$ and
that for $\alpha$,$\,\beta\in\Gamma$ the automorphisms $\alpha\circ\beta$
and $\beta\circ\alpha$ are unitarily equivalent. 
Furthermore, $\al X.\cong\Gamma/\mathrm{int}\al A.$ and
the set $\Gamma$ is sometimes called {\em dual action} on $\al A.$.
\end{rem}

For any $\gamma_1,\gamma_2\in\Gamma$ we write
\[
\gamma_1\circ\gamma_2\circ\gamma_1^{-1}\circ
\gamma_2^{-1}=\mbox{ad}\,\epsilon(\gamma_{1},\gamma_{2}),
\]
where $\epsilon(\gamma_{1},\gamma_{2})\in {\cal U}({\cal A})$ and
the class $\widehat{\epsilon}(\gamma_{1},\gamma_{2}):=
\epsilon(\gamma_{1},\gamma_{2})\,\mbox{\rm mod}\,{\cal U}({\cal Z})$ is
uniquely defined.

\begin{lem}\label{Lem3}
The permutators $\epsilon(\cdot,\cdot)$ satisfy
the following relations: 
\begin{eqnarray*} 
\epsilon(\gamma_{1},\gamma_{2})\epsilon(\gamma_{2},\gamma_{1})
&\equiv&\EINS\,\,\mbox{\rm mod}\,{\cal U}({\cal Z})\,,
        \qquad\gamma_1,\gamma_2\in\Gamma\,,\\
\epsilon(\iota,\gamma)\equiv \epsilon(\gamma,\iota)
&\equiv&\EINS\,\,\mbox{\rm mod}\,{\cal U}({\cal Z})\,,
         \qquad\gamma\in\Gamma\,,\\
\gamma_1(\epsilon(\gamma_{2},\gamma_{3}))
\epsilon(\gamma_{1},\gamma_{3})
&\equiv& \epsilon(\gamma_{1}\gamma_{2},\gamma_{3})
         \,\,\mbox{\rm mod}\,{\cal U}({\cal Z})\,,
         \qquad\gamma_1,\gamma_2,\gamma_3\in\Gamma\,,\\
A\gamma_1(B)\epsilon(\gamma_{1},\gamma_{2})
&\equiv&\epsilon(\gamma_{1}',\gamma_{2}')B\gamma_2(A)\,\,
        \mbox{\rm mod}\,{\cal U}({\cal Z})\,,\qquad
        \gamma_1,\gamma_2,\gamma_1',\gamma_2'\in\Gamma\;\;\mathrm{and}\\
&&\qquad\qquad\qquad\quad
               A\in(\gamma_1,\gamma_1')\cap\al U.(\al A.)\,,\; 
               B\in (\gamma_2,\gamma_2')\cap\al U.(\al A.).
\end{eqnarray*}
\end{lem}
\begin{beweis}
The first and second equations above are obvious. To prove the third one 
consider the the inner automorphism characterized by the l.h.s.~of 
the equation:
\begin{eqnarray*}
 \mathrm{ad}\Big(\gamma_1(\epsilon(\gamma_{2},\gamma_{3}))
                 \epsilon(\gamma_{1},\gamma_{3})\Big)
    &=& \mathrm{ad}\Big(\gamma_1(\epsilon(\gamma_{2},\gamma_{3}))\Big)\circ
        \mathrm{ad}\Big(\epsilon(\gamma_{1},\gamma_{3})\Big)\\
    &=& \gamma_1\,\mathrm{ad}(\epsilon(\gamma_{2},\gamma_{3}))\,
        \gamma_1^{-1}\circ\mathrm{ad}(\epsilon(\gamma_{1},\gamma_{3}))\\
    &=& \gamma_1\,(\gamma_2\gamma_3\gamma_2^{-1}\gamma_3^{-1})
        \,\gamma_1^{-1}\,(\gamma_1\gamma_3\gamma_1^{-1}\gamma_3^{-1})
        = (\gamma_1\gamma_2)\,\gamma_3\,(\gamma_1\gamma_2)^{-1}
           \,\gamma_3^{-1} \\
    &=& \mathrm{ad}\Big(\epsilon(\gamma_{1}\gamma_{2},\gamma_{3})\Big)\,,
\end{eqnarray*}
and this shows the desired relation. 
Finally, to prove the last equation recall that from the assumptions we have
$\gamma_1'=\mathrm{ad}(A)\circ\gamma_1$ and
$\gamma_2'=\mathrm{ad}(B)\circ\gamma_2$.
From this we compute 
\begin{eqnarray*}
 \mathrm{ad}(\epsilon(\gamma_1',\gamma_2'))
    &=& (\mathrm{ad}(A)\circ\gamma_1)\circ
        (\mathrm{ad}(B)\circ\gamma_2)\circ
        (\mathrm{ad}(A)\circ\gamma_1)^{-1}\circ
        (\mathrm{ad}(B)\circ\gamma_2)^{-1}               \\
    &=&  \mathrm{ad}(A)\circ\mathrm{ad}(\gamma_1(B))\circ
         \underbrace{\gamma_1\circ\gamma_2\circ
         \gamma_1^{-1}\circ\gamma_2^{
         -1}}_{\mathrm{ad}(\epsilon (\gamma_1,\gamma_2))}\,
         \circ\,\mathrm{ad}(\gamma_2(A))^{-1}\circ
         \mathrm{ad}(B)^{-1} \,.
\end{eqnarray*}
Therefore we get
\[
  \mathrm{ad}\Big(\epsilon(\gamma_1',\gamma_2')B\gamma_2(A)\Big)
 =\mathrm{ad}\Big(A\gamma_1(B)\epsilon (\gamma_1,\gamma_2)\Big)
\]
which implies the last equation of the statement.
\end{beweis}

\begin{defi}\label{lift}
Let $\beta_{\chi}\in\Theta(\chi),\,\chi\in{\cal
X}$, with $\beta_{\iota}=\mbox{id}_{\cal A},$ be a system of
representatives, i.e. $\pi(\beta_{\chi})=\Theta(\chi)$. Then
$\beta_{\cal X}$ is called a {\em lifting} of $\Theta$ if ${\cal
X}\ni\chi \rightarrow \beta_{\chi}\in\mbox{aut}\,{\cal A}$ is a
homomorphism.
\end{defi}

\begin{rem}
For the notion of lifting see for example Jones \cite{Jones80}.
Sutherland \cite{Sutherland80,SutherlandIn82} says that $\Theta$ {\em
splits} if there is a lifting of $\Theta$. If $\Theta$ is an isomorphism
then a lifting is also called {\em monomorphic section} (this 
latter name is used by Doplicher/Haag/Roberts \cite{DHR69b}). 

Results on the existence of liftings when ${\cal A}$ 
is a von Neumann algebra and in a more general context w.r.t.~the
group ${\cal X}$ (theory of Q-kernels) are due to Sutherland 
\cite{Sutherland80,SutherlandIn82}. Further, recall also the result of
Doplicher/Haag/Roberts \cite{DHR69b} in the ``automorphism case'' 
of the superselection theory, where ${\cal Z}=\C\,\EINS$ and 
${\cal A}$ is a so-called quasilocal algebra w.r.t. a net of 
local von Neumann algebras (see also \cite{bBaumgaertel95}).
\end{rem}

\section{Hilbert extensions}\label{Sec3}

The question concerning the description of $\{{\cal
F},\alpha_{\cal G}\}$ by ${\cal A}$ and `something else' is called
the {\em reconstruction problem}. It is posed, for example, by
Takesaki \cite[p.~202]{TakesakiIn84} and by  Bratteli/Robinson 
\cite[p.~137]{bBratteli87}. Also
the superselection structures in algebraic quantum
field theory are connected with the reconstruction problem (for
the automorphism case see Doplicher/Haag/Roberts \cite{DHR69b}).

From Lemma~\ref{Lem1} it seems natural to consider 
the corresponding {\em inverse problem}, 
which is an extension problem. This is just the
emphasis in the mentioned papers by Sutherland and Jones (see
also Nakamura/Takeda \cite{Nakamura59,Takeda59}) as well as
an essential aspect of the superselection theory
(cf.~\cite{DHR69b,bBaumgaertel95}).

\begin{defi}
Let a system $\{{\cal A},\Theta({\cal X})\}$ be
given where ${\cal X}$ is a discrete abelian group and where 
$\Theta\colon\ {\cal X}\rightarrow \mbox{Out}\,{\cal A}$ is a homomorphism
and put ${\cal G}:=\hat{\cal X}$. A Hilbert
C*-system $\{{\cal F},\alpha_{\cal G}\}$ is called a {\em
Hilbert extension} of $\{{\cal A},\Theta({\cal X})\}$ if ${\cal
A}= \Pi_{\iota}{\cal F}$ and $\Theta({\cal X})$ coincides with
the homomorphism given by Lemma~\ref{Lem1}.
\end{defi}

Now let $\{{\cal A},\Theta({\cal X})\}$ and ${\cal G}$ be given as in
the previous definition.
As it is pointed out, for example in \cite{Jones80}, a
crucial object for the extension problem is the so-called {\em
obstruction} Ob$\,\Theta$. We recall the relevant relations:
Choose a system  $\beta_{\chi}\in \Theta(\chi),\,\chi\in {\cal
X}, \beta_{\iota}:=\mbox{id}_{\cal A}$ of representatives. Then

\begin{equation}\label{ad2}
\beta_{\chi_{1}}\circ\beta_{\chi_{2}}=\mbox{ad}\,(\omega(\chi_{1},\chi_{2}))
\circ\beta_{\chi_{1}\chi_{2}},
\end{equation}  
where
\begin{equation}\label{ad3}
{\cal X}\times {\cal X}\ni (\chi_{1},\chi_{2})\rightarrow
\omega(\chi_{1},\chi_{2})\in {\cal U}({\cal A})
\end{equation}
and we have the intertwining property
\begin{equation}\label{ad4}
\omega(\chi_{1},\chi_{2})\in (\beta_{\chi_{1}\chi_{2}},
\beta_{\chi_{1}}\circ\beta_{\chi_{2}}), 
\end{equation}
which is implied by (\ref{ad2}). 
Moreover we have
\begin{equation}\label{boundary}
\omega(\iota,\chi)=\omega(\chi,\iota)=\EINS.
\end{equation}
Now associativity yields
\[
\mbox{ad}\,(\omega(\chi_{1},\chi_{2})\omega(\chi_{1}\chi_{2},\chi_{3}))=
\mbox{ad}\,(\beta_{\chi_{1}}(\omega(\chi_{2},\chi_{3}))
\omega(\chi_{1},\chi_{2}\chi_{3}))
\]
so that there is $\gamma(\chi_{1},\chi_{2},\chi_{3})\in {\cal U}
({\cal Z})$ with
\[
\omega(\chi_{1},\chi_{2})\omega(\chi_{1}\chi_{2},\chi_{3})=
\gamma(\chi_{1},\chi_{2},\chi_{3})\beta_{\chi_{1}}(\omega(\chi_{2},\chi_{3}))
\omega(\chi_{1},\chi_{2}\chi_{3}).
\]
If $\gamma(\chi_{1},\chi_{2},\chi_{3})=\EINS$ for all
$\chi_{1},\chi_{2},\chi_{3}\in{\cal X}$ we obtain the equation
\begin{equation}\label{cocycle}
\omega(\chi_{1},\chi_{2})\omega(\chi_{1}\chi_{2},\chi_{3})=
\beta_{\chi_{1}}(\omega(\chi_{2},\chi_{3}))\omega(\chi_{1},\chi_{2}\chi_{3}).
\end{equation}

Obviously, the existence of a system of representatives
$\beta_{\cal X}$ such that equation (\ref{cocycle}) has a solution $\omega$
equipped with the properties (\ref{ad3})--(\ref{boundary}) is necessary for the
existence of a Hilbert extension.
Even more, the existence of such a solution is also sufficient for
the existence of a Hilbert extension.

\begin{defi}
A function $\omega$, assigned to a given system
$\beta_{\cal X}$ of representatives of $\Theta({\cal X})$,
equipped with the properties (\ref{ad3})--(\ref{cocycle}) is called a generalized 2-cocycle.
\end{defi}

One calculates easily that the existence of a generalized
2-cocycle is independent of the choice of the system
$\beta_{\cal X}$ of representatives. Further, a generalized cocycle
$\omega$ for $\beta_{\cal X}$ satisfies the relation
\[
\mbox{ad}\,(\omega(\chi_{1},\chi_{2})\omega(\chi_{2},\chi_{1})^{-1})=
\beta_{\chi_{1}}\circ\beta_{\chi_{2}}\circ\beta_{\chi_{1}}^{-1}\circ
\beta_{\chi_{2}}^{-1}.
\]

The existence of a lifting of $\Theta$ can be expressed in terms
of generalized 2-cocycles as follows.

\begin{lem}\label{Lem4}
There exists a lifting $\beta_{\cal X}$ of
$\Theta$ iff to each system $\gamma_{\cal X}$  of
representatives there corresponds a generalized 2--cocycle
$\omega$ of the form
\[
\omega(\chi_{1},\chi_{2})\equiv\gamma_{\chi_{1}}(V_{\chi_{2}}^{-1})
V_{\chi_{1}}^{-1}V_{\chi_{1}\chi_{2}}\quad\mathrm{mod}\,\al U.(\al Z.)\,,
\]
where $V_{\chi}\in {\cal U}({\cal A}),\,V_\iota=\EINS.$
In this case, i.e. if there is a lifting $\beta_{\cal X}$,
then a corresponding generalized 2--cocycle $\omega$ 
is given by $\omega(\chi_{1},\chi_{2})=\EINS$ for all 
$\chi_{1},\chi_{2}\in{\cal X}.$
\end{lem}
\begin{beweis}
Let $\beta_\chi= \mathrm{ad}(V_\chi)\circ\gamma_\chi$, $V_\chi\in
\al U.(\al A.)$, $\chi\in\al X.$. Now if 
$\omega(\chi_{1},\chi_{2})=\gamma_{\chi_{1}}(V_{\chi_{2}}^{-1})
V_{\chi_{1}}^{-1}V_{\chi_{1}\chi_{2}}Z$
for some $Z\in\al U.(\al Z.)$, then we have on the one hand
$\beta_{\chi_1\chi_2}= \mathrm{ad}(V_{\chi_1\chi_2})\circ
\gamma_{\chi_1\chi_2}$ and on the other
\[
   \beta_{\chi_1}\circ\beta_{\chi_2} 
  =(\mathrm{ad}(V_{\chi_1})\circ\gamma_{\chi_1})\circ
    (\mathrm{ad}(V_{\chi_2})\circ\gamma_{\chi_2})
  = \mathrm{ad}\Big(V_{\chi_1}\gamma_{\chi_1}(V_{\chi_2})\,
    \omega(\chi_1,\chi_2)\Big)\circ\gamma_{\chi_1\chi_2}\,,
\]
which using the assumption on $\omega$  and the fact that 
$\mathrm{ad}(V_{\chi_1\chi_2}Z)=\mathrm{ad}(V_{\chi_1\chi_2})$, implies that
$\beta_{\chi_1\chi_2}=\beta_{\chi_1}\circ\beta_{\chi_2}$, i.e.~there
is a lift of $\Theta$. To prove the converse let
$\beta_{\chi_1\chi_2}=\beta_{\chi_1}\circ\beta_{\chi_2}$, so that from
the above relations we have
\[
 \mathrm{ad}(V_{\chi_1\chi_2})
  = \mathrm{ad}\Big(V_{\chi_1}\gamma_{\chi_1}(V_{\chi_2})\,
    \omega(\chi_1,\chi_2)\Big)\,,
\]
which implies $\omega(\chi_{1},\chi_{2})=\gamma_{\chi_{1}}(V_{\chi_{2}}^{-1})
V_{\chi_{1}}^{-1}V_{\chi_{1}\chi_{2}}\;\mathrm{mod}\,\al U.(\al Z.)$.
\end{beweis}

\begin{teo}\label{Main1}
Let $\omega$ be a generalized 2--cocycle
for the system $\beta_{\cal X}$ of representatives. Then
there is a Hilbert extension $\{{\cal F},\alpha_{\cal G}\}$
of $\{{\cal A},\Theta({\cal X})\}.$
\end{teo}
\begin{beweis}
The proof consists of several steps that correspond to gradually
imposing a richer structure on an initially considered $\al A.$--left
module:

1.~Indeed, choose first
system of 1-dimensional linear spaces, generated by abstract
elements $U_{\chi}$, $\chi\in\al X.$, $U_{\iota}:=\EINS\in {\cal A}$. 
Form the ${\cal A}$--left modules ${\cal A}\otimes\C U_{\chi}$ and
${\cal F}_{0}:=\bigoplus_{\chi}({\cal
A}\otimes\C U_{\chi}).$ 
By identification $A\otimes\EINS\leftrightarrow A,\,\EINS\otimes
U_{\chi}\leftrightarrow U_{\chi}$ one has
\[
{\cal F}_ {0}=\Bigg\{\sum_{\mbox{\tiny $\chi\;$, finite~sum}}
              A_{\chi}U_{\chi}\mid \quad A_{\chi}\in{\cal A}\Bigg\}\,,
\]
where $\{ U_\chi\mid\chi\in\al X.\}$ forms an 
abstract $\al A.$--module basis. 

2.~Next we want to equip $\al F._0$ with a multiplication structure.
First $\al F._0$ becomes an $\al A.$--bimodule extending linearly
the following definition
\[
 U_\chi\,A:=\beta_\chi(A)\,U_\chi\,,\quad A\in\al A.\,,\chi\in\al X.\,,
\]
where $\beta_\al X.$ is the system of representatives to which we
associate the generalized cocycle $\omega$. Now the product structure
is finally specified by putting 
\[
 U_{\chi_1}\cdot U_{\chi_2}:=\omega(\chi_1,\chi_2)\,U_{\chi_1\chi_2}\,,
 \quad \chi_1,\chi_2\in\al X.\,,
\]
where the cocycle equation (\ref{cocycle}) guarantees that the product
is associative and the boundary conditions (\ref{boundary}) lead to
$U_\chi\cdot\EINS=\EINS\cdot U_\chi=U_\chi$. Note that the preceding 
product structure already implies that the $U_\chi$ are invertible.
Indeed, it can be checked easily that the inverse is given explicitly
by 
\[
 U_\chi^{-1}:=\beta_{\chi^{-1}}\Big(\omega(\chi,\chi^{-1})^{-1}\Big)
              U_{\chi^{-1}}
\]
(use for example the relation $\beta_\chi(\omega(\chi^{-1},\chi))=
\omega(\chi,\chi^{-1})$, which follows from the cocycle equation 
(\ref{cocycle}) by putting $\chi_1:=\chi$, $\chi_2:=\chi^{-1}$
and $\chi_3=\chi$).

3.~The following step consists in defining a *--structure on 
$\al F._0$. This is done by putting  
\[
 U_\chi^*:=\omega(\chi^{-1},\chi)^* U_{\chi^{-1}}\quad\mathrm{and}
 \quad (AU_\chi)^*:=U_\chi^*A^*.
\]
We still have to check that this definition is consistent,
in particular with the product structure in $\al F._0$, i.e.~we
have to verify:
\begin{equation}\label{consistent}
 (U_\chi^*)^*=U_\chi\,,\quad (U_\chi A)^*=A^*U_\chi^*\quad
 \mathrm{and}\quad (U_{\chi_1}\cdot U_{\chi_2})^*=
 U_{\chi_2}^*\cdot U_{\chi_1}^*\,.
\end{equation}
For the first equation we have
\begin{eqnarray*}
 (U_\chi^*)^* &=& \Big(\omega(\chi^{-1},\chi)^*\, U_{\chi^{-1}}\Big)^*
               = U_{\chi^{-1}}^*\,\omega(\chi^{-1},\chi)
               =\omega(\chi,\chi^{-1})^* \,U_\chi\,\omega(\chi^{-1},\chi)\\
              &=& \omega(\chi,\chi^{-1})^*\,\beta_\chi\Big(
                  \omega(\chi^{-1},\chi)^*\Big) U_\chi
               = \omega(\chi,\chi^{-1})^*\,\omega(\chi,\chi^{-1})\, U_\chi\\
              &=& U_\chi
\end{eqnarray*}
The second equation in (\ref{consistent}) can also be checked immediately
from the definitions considered above. For the last equation we will
consider the two sides separately: for the r.h.s.~we have
\begin{eqnarray*}
 U_{\chi_2}^*\cdot U_{\chi_1}^*
        &=& \omega(\chi_2^{-1},\chi_2)^*\,U_{\chi_2^{-1}}\cdot
             \omega(\chi_1^{-1},\chi_1)^*\,U_{\chi_1^{-1}}\\
        &=& \omega(\chi_2^{-1},\chi_2)^*\,\beta_{\chi_2^{-1}}\Big(         
          \omega(\chi_1^{-1},\chi_1)^*\Big)U_{\chi_2^{-1}} \,U_{\chi_1^{-1}}\\
        &=& \omega(\chi_2^{-1},\chi_2)^*\,\beta_{\chi_2^{-1}}\Big(         
          \omega(\chi_1^{-1},\chi_1)^*\Big)\omega(\chi_2^{-1},\chi_1^{-1})
          \,U_{(\chi_1\chi_2)^{-1}}\\
        &=& \omega(\chi_2^{-1},\chi_2)^*\,\omega((\chi_1\chi_2)^{-1},\chi_1)^*
          \underbrace{\omega(\chi_2^{-1},\chi_1^{-1})^*
         \omega(\chi_2^{-1},\chi_1^{-1})}_{\EINS}\,U_{(\chi_1\chi_2)^{-1}}\,,
\end{eqnarray*}
where we have used the relation
\[
 \beta_{\chi_2^{-1}}(\omega(\chi_1^{-1},\chi_1))=
\omega(\chi_2^{-1},\chi_1^{-1})\,\omega(\chi_2^{-1}\chi_1^{-1},\chi_1)\,,
\]
which again follows from the cocycle equation 
(\ref{cocycle}) taking now $\chi_1:=\chi_2^{-1}$, $\chi_2:=\chi_1^{-1}$
and $\chi_3=\chi_1$. Now the l.h.s.~reads
\begin{eqnarray*}
 (U_{\chi_1}\cdot U_{\chi_2})^*
        &=& U_{\chi_1\chi_2}^*\,\omega(\chi_1,\chi_2)^*
      =\omega((\chi_1\chi_2)^{-1},\chi_1\chi_2)^*\,\,U_{(\chi_1\chi_2)^{-1}}
       \,\omega(\chi_1,\chi_2)^*\\
        &=& \omega((\chi_1\chi_2)^{-1},\chi_1\chi_2)^*
            \beta_{(\chi_1\chi_2)^{-1}}\Big(\omega(\chi_1,\chi_2)^*\Big)
            \,U_{(\chi_1\chi_2)^{-1}}\,.
\end{eqnarray*}
Thus to show the last equation in (\ref{consistent}) we need to prove
that 
\[
 \omega((\chi_1\chi_2)^{-1},\chi_1\chi_2)^*\,
 \beta_{(\chi_1\chi_2)^{-1}}\Big(\omega(\chi_1,\chi_2)^*\Big)
  =\omega(\chi_2^{-1},\chi_2)^*\,
    \omega\Big((\chi_1\chi_2)^{-1},\chi_1\Big)^*
\]
or taking adjoints
\[
 \beta_{(\chi_1\chi_2)^{-1}}\Big(\omega(\chi_1,\chi_2)\Big)
 \,\omega((\chi_1\chi_2)^{-1},\chi_1\chi_2)
 = \omega\Big((\chi_1\chi_2)^{-1},\chi_1\Big)\,
   \omega(\chi_2^{-1},\chi_2)\,.
\]
But the preceding equation is nothing else than the cocycle equation
(\ref{consistent}) with $\chi_1:=(\chi_1\chi_2)^{-1}$, 
$\chi_2:=\chi_1$ and $\chi_3:=\chi_2$. Finally, note that since
$\beta_{\chi^{-1}}\Big(\omega(\chi,\chi^{-1})^{-1}\Big)
=\omega(\chi^{-1},\chi)^*$ we also have that the $U_\chi$, 
are unitary, i.e.~$U_\chi^*=U_\chi^{-1}$, $\chi\in\al X.$.

4.~Here we will define a representation of the compact abelian group 
$\al G.=\widehat{\al X.}$ in terms of automorphisms of the *--algebra
$\al F._0$. The automorphisms are fixed by putting
\[
 \alpha_g(U_\chi):=\chi(g)\,U_\chi\quad\mathrm{and}\quad
 \alpha_g(AU_\chi):=A\,\alpha_g(U_\chi)=\chi(g)\, A\,U_\chi\,,\quad
   g\in\al G.,\,A\in\al A.,\,\chi\in\al X.\,.
\]
First we check that with the definition above the $\alpha_g$ is indeed 
an automorphism compatible with the structure in $\al F._0$:
\begin{eqnarray*}
\alpha_g\Big(U_{\chi_1}U_{\chi_2}\Big)
 &=& \alpha_g\Big(\omega(\chi_1,\chi_2)\,U_{\chi_1\chi_2}\Big)
     =(\chi_1\chi_2)(g)\,\omega(\chi_1,\chi_2)\,U_{\chi_1\chi_2}\\
 &=& \chi_1(g)\chi_2(g)\,U_{\chi_1}\,U_{\chi_2}
     =\alpha_g\Big(U_{\chi_1}\Big)\alpha_g\Big( U_{\chi_2}\Big)
\end{eqnarray*}
and
\begin{eqnarray*}
\alpha_g\Big(U_\chi^*\Big)
 &=& \alpha_g\Big(\omega(\chi^{-1},\chi)^*\,U_{\chi^{-1}}\Big)
     =(\chi^{-1})(g)\,\omega(\chi^{-1},\chi)^*\,U_{\chi^{-1}}\\
 &=& \overline{\chi}(g) \,U_\chi^*
     =\alpha_g\Big(U_\chi\Big)^*\,.
\end{eqnarray*}
It can be also easily seen that the assignment $\al G.\ni g\to
\alpha_g\in\mathrm{aut}\,\al F._0$ is an injective group homomorphism.
Finally, note that the fixed point algebra of the previous action
coincides with $\al A.$, i.e.~for $F\in\al F._0$, $\alpha_g(F)=F$
for all $g\in\al G.$ iff $F\in\al A.$. Indeed, for an arbitrary
element $\sum_{\chi}A_{\chi}U_{\chi}\in\al F._0$ the equation
$\sum_{\chi}\chi(g)A_{\chi}U_{\chi}=\sum_{\chi}A_{\chi}U_{\chi}$,
$g\in\al G.$, implies by the base property of the $U_\chi$ that
$\chi(g)A_{\chi}=A_{\chi}$ , $g\in\al G.$, $\chi\in\al X.$. Therefore
if $\chi_0\not=\iota$, then there is a $g_0\in\al G.$ with $\chi_0(g_0)
\not= 1$ and this shows that $A_{\chi_0}=0$. The converse implication
is obvious.

5.~Finally, to specify a C*--norm on $\al F._0$ we introduce
the following ${\cal A}$--valued scalar product (note the variation
w.r.t.~the definition in \cite[p.~101]{bBaumgaertel95}):
\[
\langle F_{1},F_{2}\rangle:=
\sum_{\chi}\beta_{\chi}^{-1}(A_{\chi}^{\ast}B_{\chi})\,,\quad\mathrm{where}
\quad F_{1}=\sum_{\chi}A_{\chi}U_{\chi}\,,\;F_{2}=\sum_{\chi}B_{\chi}U_{\chi} 
\in\al F._0. 
\]
This scalar product satisfies the properties
\[
 \langle F_{1},F_{2}\rangle^*=\langle F_{2},F_{1}\rangle\,,\quad
 \langle F_1,F_1\rangle\geq 0\quad\mathrm{and}\quad
 \langle F_1,F_1\rangle=0\;\mathrm{iff}\;F_1=0\,.
\]
Next we show that 
\[
 \langle F_{1},F_{2}\rangle=\Pi_\iota(F_1^*F_2)\,,
\]
Indeed, using the definitions above we have
\[
 F_1^*F_2=\sum_{\chi_1,\chi_2}U_{\chi_1}^* A_{\chi_1}^*B_{\chi_2} 
          U_{\chi_2}
         =\sum_{\chi_1,\chi_2}\omega(\chi_1^{-1},\chi_1)^*\,
          \beta_{\chi_1^{-1}}\big(A_{\chi_1}^*B_{\chi_2}\big)\,
          \omega(\chi_1^{-1},\chi_2)U_{\chi_1^{-1}\chi_2}\,.
\]
Putting, $\chi_1=\chi_2=\chi$ in the preceding expression we get
\begin{eqnarray*}
\Pi_\iota(F_1^*F_2)
  &=& \sum_{\chi}\omega(\chi^{-1},\chi)^*\,\beta_{\chi^{-1}}
      \big(A_{\chi}^*B_{\chi}\big)\,\omega(\chi^{-1},\chi)\\ 
  &=& \sum_{\chi}\omega(\chi^{-1},\chi)^*\omega(\chi^{-1},\chi)
      \,\beta_{\chi}^{-1}\big(A_{\chi}^*B_{\chi}\big)
      \,\omega(\chi^{-1},\chi)^*\omega(\chi^{-1},\chi)\\
  &=& \langle F_{1},F_{2}\rangle\,,
\end{eqnarray*}
where for the second equation before we have used eq.~(\ref{ad2})
in the form $\beta_{\chi^{-1}}=\mathrm{ad}\,(\omega(\chi^{-1},\chi))
\circ\beta_{\chi}^{-1}$. In particular the relation above implies
the  following invariance property: $\langle \alpha_g(F_{1}),
\alpha_g(F_{2})\rangle=\langle F_{1},F_{2}\rangle$, $g\in\al G.$.

Define next the following norm on $\al F._0$ by 
\[
 |F|:=\|\langle F,F\rangle\|^{\frac12}\,,\quad F\in\al F._0\,,
\]
and the representation of $\al F._0$ on $(\al F._0,|\cdot|)$ in terms of
multiplication operators
\[
 \rho(F)X:=FX\,,\quad F,X\in\al F._0\,.
\]
Note that by the definition of the $\al A.$--valued scalar product
the property $\rho(F^*)=\rho(F)^*$, $F\in\al F._0$, holds. Now
using the corresponding operator norm we introduce
\[
\|F\|_*:=|\rho(F)|_{op}\,,\quad F\in\al F._0\,,
\]
which by similar arguments as in \cite[p.~102-103]{bBaumgaertel95}
satisfies the C*--property $\|F^*F\|_*=\|F\|_*^2$. Further, it satisfies 
also (cf.~again the previous reference)
\[
 \|A\|_*=\|A\|\,,\quad A\in\al A.\quad\mathrm{and}\quad
 \|\alpha_g(F)\|_*=\|F\|_*\,,\quad g\in\al G.\,,\, F\in\al F._0\,.
\]
Therefore, we can finally extend $\alpha_g$ isometrically from $\al F._0$
to
\[
 \al F.:=\mathrm{clo}_{\|\cdot\|_*}(\al F._0)\,.
\] 
Further, $\alpha_\al G.\subset\mathrm{aut}\,\al F.$ is norm continuous
w.r.t.~the pointwise norm convergence, because for any 
$F_0=\sum_\chi A_\chi U_\chi\in\al F._0$ we have
\[
  \|\alpha_{g_1}(F_0)-\alpha_{g_2}(F_0)\|_*
    = \|\sum_\chi \Big(\chi(g_1)-\chi(g_2)\Big)\,A_\chi U_\chi\|_*
    \leq \sum_\chi |\chi(g_1)-\chi(g_2)|\,\|A_\chi\|\,.
\]
By construction we also have that $U_\chi\in\Pi_\chi(\al F.)$, $\chi\in\al
X.$. Therefore from the definitions of Sections~\ref{Sec2} and
\ref{Sec3} we have constructed a Hilbert C*--extension
$\{\al F.,\alpha_\al G.\}$ of $\{\al A.,\Gamma\}$ and the proof is
concluded.
\end{beweis}

Using now Lemma~\ref{Lem4} one has
\begin{cor}\label{1lift}
If there is a lifting of $\Theta$, then
there is a Hilbert extension of $\{{\cal A},\Theta({\cal X})\}$,
corresponding to $\omega =\EINS$.
\end{cor}
\begin{rem}
The construction in the proof of the previous theorem generalizes
to the nontrivial center situation the procedure already presented
(with small modifications) in \cite[Section~3.6]{bBaumgaertel95}.
\end{rem}

The second problem consists in the description of all Hilbert
extensions. For this purpose let 
$\Omega({\cal X},{\cal U}({\cal Z}),\beta_{\cal X})$
be the set of all ${\cal U}({\cal Z})$--valued 2--cocycles
$\lambda$, i.e. $\lambda$ satisfies equation (\ref{cocycle}) and condition
(\ref{boundary}), but (\ref{ad3}),(\ref{ad4}) are replaced by $\lambda(\chi_{1},\chi_{2})\in 
{\cal U}({\cal Z})$. For example,
$\lambda(\chi_{1},\chi_{2}):=\EINS$ for all
$\chi_{1},\chi_{2}\in{\cal X}$
is such a cocycle.
Further let $\Omega_{0}({\cal X},{\cal U}({\cal Z}),\beta_{\cal X})$ be the
set of all ${\cal U}({\cal Z})$--valued coboundaries $\partial Z$, i.e.
\[
\partial Z(\chi_{1},\chi_{2}):=
\frac{Z(\chi_{1})\beta_{\chi_{1}}(Z(\chi_{2}))}{Z(\chi_{1}\chi_{2})}\,,
\]
where $Z(\cdot)$ is a ${\cal U}({\cal Z})$-valued 1-cycle, $Z(\iota)=\EINS$.
Then $\partial Z$ is a ${\cal
U}({\cal Z})$-valued 2-cocycle, $\Omega\supseteq \Omega_{0}$.
As usual, $\Omega$ and $\Omega_{0}$ are abelian groups w.r.t.~pointwise 
multiplication and the second cohomology is given by
$H^{2}({\cal X},{\cal U}({\cal Z}),\beta_{\cal X}):=\Omega/\Omega_{0}.$

Next we need the concept of ${\cal A}$--{\em module isomorphism}
of Hilbert extensions.

\begin{defi}
Let $\{{\cal F}^{1},\alpha_{\cal G}^{1}\},\,
\{{\cal F}^{2},\alpha_{\cal G}^{2}\}$
be Hilbert extensions of
$\{{\cal A},\Theta({\cal X})\}$.
They are called ${\cal A}$--module isomorphic if there is an
algebraic isomorphism $\Phi\colon\ \al F.^1 \to \al F.^2$,
with $\Phi(A)=A$ for all $A\in{\cal A}$ and
$\Phi\circ\alpha_{g}^{1}=\alpha_{g}^{2}\circ\Phi$ for all
$g\in{\cal G}$.
\end{defi}
\begin{teo}\label{Teo2}
Let $\omega_{0}$ be a generalized 2--cocycle. Then:
\begin{itemize}
\item[(i)] Each ${\cal U}({\cal Z})$--valued 2--cocycle
$\lambda$ yields a Hilbert extension generated by the
generalized 2--cocycle $\omega:=\lambda\cdot\omega_{0}$
and each Hilbert extension is generated by some
${\cal U}({\cal Z})$--valued 2--cocycle $\lambda$ via
$\omega:=\lambda\cdot\omega_{0}.$
\item[(ii)] Two Hilbert extensions are ${\cal A}$--module
isomorphic iff the generating generalized 2--cocycles
$\omega_{1},\omega_{2}$ differ only by a
${\cal U}({\cal Z})$--valued coboundary $\partial Z$, i.e.
$\omega_1=\partial Z\cdot\omega_2.$
\end{itemize}
\end{teo}
\begin{beweis}
(i) If two
generalized 2--cocycles $\omega_{1},\omega_{2}$ are given, then 
note first that $\lambda(\chi_{1},\chi_{2}):=
\omega_{1}(\chi_{1},\chi_{2})\omega_{2}(\chi_{1},\chi_{2})^{-1}
\in {\cal U}({\cal Z})$ for all $\chi_{1},\chi_{2}$,
because of condition (\ref{ad4}). Further, eq.~(\ref{boundary})
follows from the corresponding properties of $\omega_1$ and 
$\omega_2$. Finally, the cocycle equation for 
$\lambda(\chi_{1},\chi_{2})$ is a consequence of the following 
computation:
\begin{eqnarray*}
\lambda(\chi_{1},\chi_{2})\lambda(\chi_{1}\chi_{2},\chi_{3})
&=& \omega_{1}(\chi_{1},\chi_{2})\omega_{2}(\chi_{1},\chi_{2})^{-1}\cdot
   \omega_{1}(\chi_{1}\chi_{2},\chi_{3})\omega_{2}(\chi_{1}\chi_{2},
   \chi_{3})^{-1} \\
&=& \omega_{1}(\chi_{1},\chi_{2})\omega_{1}(\chi_{1}\chi_{2},\chi_{3})
   \omega_{2}(\chi_{1}\chi_{2},\chi_{3})^{-1}
   \omega_{2}(\chi_{1},\chi_{2})^{-1} \\
&=& (\omega_{1}(\chi_{1},\chi_{2})\omega_{1}(\chi_{1}\chi_{2},\chi_{3}))\cdot
  (\omega_{2}(\chi_{1},\chi_{2})\omega_{2}(\chi_{1}\chi_{2},\chi_{3}))^{-1}\\ 
&=& \beta_{\chi_{1}}(\omega_{1}(\chi_{2},\chi_{3}))\omega_{1}(\chi_{1},
    \chi_{2}\chi_{3})\cdot
    (\beta_{\chi_{1}}(\omega_{2}(\chi_{2},\chi_{3}))\omega_{2}(\chi_{1},
    \chi_{2}\chi_{3}))^{-1}\\
&=& \beta_{\chi_{1}}(\omega_{1}(\chi_{2},\chi_{3}))\omega_{1}(\chi_{1},
    \chi_{2}\chi_{3})
    \omega_{2}(\chi_{1},\chi_{2}\chi_{3})^{-1}\beta_{\chi_{1}}
    (\omega_{2}(\chi_{2},\chi_{3}))^{-1}\\
&=& \beta_{\chi_{1}}(\omega_{1}(\chi_{2},\chi_{3})\omega_{2}(\chi_{2},
    \chi_{3})^{-1})
    \omega_{1}(\chi_{1},\chi_{2}\chi_{3})\omega_{2}(\chi_{1},
    \chi_{2}\chi_{3})^{-1}\\
&=& \beta_{\chi_{1}}(\lambda(\chi_{2},\chi_{3}))\cdot \lambda
    (\chi_{1},\chi_{2}\chi_{3}),
\end{eqnarray*}
i.e.~if one fixes a generalized 2-cocycle $\omega_{0}$, then
$\omega:=\lambda\cdot \omega_{0}$ runs through all generalized
2--cocycles $\omega$ if $\lambda$ runs through all ${\cal
U}({\cal Z})$--valued 2--cocycles in 
$\Omega({\cal X},{\cal U}({\cal Z}),\beta_{\cal X})$. 

(ii) Let $\{{\cal F}^{1},\alpha_{\cal G}^{1}\}$ and
$\{{\cal F}^{2},\alpha_{\cal G}^{2}\}$ be two Hilbert extensions of 
$\{{\cal A},\Theta({\cal X})\}$ and denote the corresponding
set of abstract unitaries by $\{U_\chi\mid \chi\in\al X.\}$
resp.~$\{V_\chi\mid \chi\in\al X.\}$.

Suppose first that there exists coboundary $\partial Z\in
\Omega_{0}({\cal X},{\cal U}({\cal Z}),\beta_{\cal X})$, where
$\beta_\al X.$ is system of representatives in $\Theta$, such
that the corresponding generalized cocycles $\omega_1$ and
$\omega_2$ satisfy $\omega_1=\partial Z\cdot \omega_2$.
In this case we will show that the extensions are isomorphic. Indeed,
define the isomorphism by
\[
 \Phi(AU_\chi):=A\,Z(\chi)\,V_\chi\,,\quad A\in\al A.\,,\;\chi\in\al X.\,,
\]
and extend it by linearity to the corresponding left $\al A.$--module.
Now $\Phi$ is even a *--homomorphism between the *--algebras
$\al F.^1_0$ and $\al F.^2_0$ that are defined in step~3 of the proof of
Theorem~\ref{Main1}. This follows from the following computations:
\begin{eqnarray*}
\Phi(U_\chi A)&=&\Phi\Big(\beta_\chi(A)U_\chi\Big)=Z(\chi)\,V_\chi A
                 =\Phi(U_\chi)\Phi(A)\,,\\[3mm]
\Phi(U_\chi U_{\chi'})&=&\Phi\Big(\omega_1(\chi,\chi')U_{\chi\chi'}\Big)
     = \partial Z(\chi,\chi')\cdot\omega_2(\chi,\chi')\,
       Z(\chi\chi')\,V_{\chi\chi'}\\
  &=& \frac{Z(\chi)\beta_{\chi}(Z(\chi'))}{Z(\chi\chi')}\cdot
       Z(\chi\chi')\,V_{\chi}V_{\chi'}=Z(\chi)V_{\chi}\,Z(\chi')V_{\chi'}
      =\Phi(U_\chi)\Phi(U_{\chi'})\,,\\[3mm]
\Phi(U_\chi^*)&=&\Phi\Big(\omega_1(\chi^{-1},\chi)^*U_{\chi^{-1}}\Big)
      = \partial Z(\chi^{-1},\chi)^*\cdot\omega_2(\chi^{-1},\chi)^*
       Z(\chi^{-1})V_{\chi^{-1}}\\
    &=& Z(\chi^{-1})^*\beta_{\chi^{-1}}(Z(\chi))^*Z(\chi^{-1})\,
        \omega_2(\chi^{-1},\chi)^*V_{\chi^{-1}}=(Z(\chi)V_\chi)^*
        =\Phi(U_\chi)^*\,,
\end{eqnarray*} 
where $\chi,\chi'\in\al X.$, $A\in\al A.$. Note further that on 
$\al F._0^1$ we already have $\Phi\circ\alpha_g^1=\alpha_g^2\circ\Phi$,
$g\in\al G.$, since for any $\chi\in\al X.$ we have
\[
 \Phi\circ\alpha_g^1(AU_\chi)=\chi(g)\, A\,Z(\chi)V_\chi
 =\alpha_g^2(A\,Z(\chi)V_\chi)=\alpha_g^2\circ\Phi(AU_\chi)\,.
\]
Recall that $\Phi$ is a bijection between $\al F._0^1$ and
$\al F._0^2$ and we will finish this part of the proof if we can
also show that $\Phi$ is even an isometry w.r.t~the corresponding
C*--norms, because in this case we can isometrically extend $\Phi$
to the desired Hilbert extension isomorphism $\Phi\colon\ \al F.^1
\to \al F.^2$. Now denote by $\langle\cdot,\cdot\rangle_k$ the
$\al A.$--valued scalar products on $\al F._0^k$, $k=1,2$, 
given in step~5 of
the proof of Theorem~\ref{Main1}. For any $F=\sum_\chi A_\chi U_\chi\in
\al F._0^1$, so that $\Phi(F)=\sum_\chi A_\chi\,Z(\chi)\, V_\chi
\in\al F._0^2$, we have the following invariance
\[
 \langle\Phi(F),\Phi(F)\rangle_2=\sum_\chi\beta_\chi^{-1}\Big(
   Z(\chi)^*A_\chi^*A_\chi Z(\chi)\Big)
 =\sum_\chi\beta_\chi^{-1}\Big(A_\chi^*A_\chi\Big)
 = \langle F,F\rangle_1\,.
\]
From this and recalling the definition of the C*--norm again in 
step~5 of the proof of Theorem~\ref{Main1} we immediately get the desired 
isometry property:
\[
\|\Phi(F)\|_*=\mathop{\mathrm{sup}}\limits_{\mbox{\tiny         
                                  $\begin{array}{c} X_2\in\al F._0^2\\[.5mm]
                                  |X_2|\leq 1\end{array}$}} |\Phi(F) X_2|
             =\mathop{\mathrm{sup}}\limits_{\mbox{\tiny 
                       $\begin{array}{c}
                        X_1\in\al F._0^1\\[.5mm]
                        |X_1|\leq 1\end{array}$}}|\Phi(F) X_1|=\|F\|_*\,.
\]

To prove the converse implication assume that $\Phi\colon\ \al F._1
\to \al F._2$ specifies the isomorphy of the Hilbert extensions.
Use the unitaries $\{U_\chi\mid\chi\in\al X. \}$ and 
$\{V_\chi\mid\chi\in\al X. \}$ in $\al F._1$ resp.~$\al F._2$ to define
the unitary
\[
 Z(\chi):=\Phi(U_\chi)\,V_\chi^*\,,\quad\chi\in\al X.\,,
\]
that satisfies $Z(\iota)=\EINS$. Even more $Z(\chi)\in\al U.(\al Z.)$,
since for any $A\in\al A.$ we have
\[
 A\,Z(\chi)=\Phi(AU_\chi)\,V_\chi^*=\Phi\Big(U_\chi\beta_{\chi}^{-1}(A)\Big)
            V_\chi^*=\Phi(U_\chi)\,(A^*V_\chi)^*=Z(\chi)\,A\,.
\]
Finally, for $\chi,\chi'\in\al X.$ we have
\begin{eqnarray*}
Z(\chi\chi')&=&\Phi\Big(\omega_1(\chi,\chi')^{-1} U_\chi U_{\chi'}\Big)
              \cdot V_{\chi'}^*V_{\chi}^*\,(\omega_2(\chi,\chi')^{-1})^*\\
            &=&\omega_1(\chi,\chi')^{-1}\,\Phi(U_\chi)\,Z(\chi')\,V_\chi^*
               \,\omega_2(\chi,\chi')\\
             &=&\omega_1(\chi,\chi')^{-1}\,\Phi(U_\chi)\Big(
              \beta_{\chi}(Z(\chi')^*)V_\chi\Big)^*\omega_2(\chi,\chi')\\
            &=& \omega_1(\chi,\chi')^{-1}\;Z(\chi)\,\beta_{\chi}(Z(\chi'))
                \;\omega_2(\chi,\chi')\,.
\end{eqnarray*}
Now recalling the definition of the coboundary $\partial Z$, the preceding
equations imply that $\omega_1(\chi,\chi')=\partial Z(\chi,\chi')
\cdot \omega_2(\chi,\chi')$, $\chi,\chi'\in\al X.$, and the prove is 
concluded.
\end{beweis}

\begin{rem}
\begin{itemize}
\item[(i)] Note that the results are independent of the choice of the system
$\beta_{\cal X}$ of representatives of $\Theta({\cal X})$.
Theorem~\ref{Teo2} means that there is a bijection between
$H^{2}({\cal X},{\cal U}({\cal Z}),\beta_{\cal X})$
and the set of all ${\cal A}$--module isomorphy classes of Hilbert
extensions of
$\{{\cal A},\Theta({\cal X})\}$ if there is one extension.
In other words, the theorem gives an {\em outer}
characterization of
$H^{2}({\cal X},{\cal U}({\cal Z}),\beta_{\cal X})$
by the set of all ${\cal A}$--module isomorphy classes of Hilbert
extensions.
\item[(ii)] For a closer analysis of the second cohomology in the 
special cases were $\Gamma\cong\Z_{N}$ and $\Gamma\cong\Z_2\times
\Z_2$ see \cite{p383}. Consider also the abstract results in
\cite[Chapter~4]{bMacLane95}.
\end{itemize}
\end{rem}

\section{The case of a trivial center}  

In this case we have ${\cal Z}=\C\EINS$, thus
${\cal U}({\cal Z})=\T\EINS$ and this implies that
two automorphisms $\alpha,\beta\in\Gamma$ are 
either unitarily equivalent or otherwise disjoint. 
The following result is a special case of the
famous Doplicher/Roberts theorem (see \cite{Doplicher89b,Baumgaertel97})
in the present automorphism context.

\begin{pro}
If there is a system of representatives $\epsilon(\alpha,\beta)$
of the permutator classes $\widehat{\epsilon}(\alpha,\beta)$ which
satisfy the equations 
\begin{eqnarray*} 
\epsilon(\gamma_1,\gamma_2)\epsilon(\gamma_2,\gamma_1)
     &=&\EINS\,,\\
\epsilon(\iota,\gamma)= \epsilon(\gamma,\iota)
     &=&\EINS\,,\\
\gamma_1(\epsilon(\gamma_2,\gamma_3))\epsilon(\gamma_1,\gamma_3)
     &=&\epsilon(\gamma_1\gamma_2,\gamma_3)\,,\\
A\beta_{\chi_{1}}(B)\epsilon(\chi_{1},\chi_{2})
&=&\epsilon'(\chi_{1},\chi_{2})B\beta_{\chi_{2}}(A)\,,
\end{eqnarray*}
for all $A\in(\beta_{\chi_{1}},\beta'_{\chi_{1}})$, 
$B\in(\beta_{\chi_{2}},\beta'_{\chi_{2}})$,
where $\epsilon'$ belongs to $\beta'_{\cal X}$, then there is a generalized
2--cocycle $\omega_{0}$ w.r.t.~some system
$\beta_\chi$ of representatives of the classes 
$\chi\in\Gamma/\mathrm{int}\al A.$, with
\[
\omega_{0}(\chi_{1},\chi_{2})\omega_{0}(\chi_{2},\chi_{1})^{-1}=
\epsilon(\beta_{\chi_{1}},\beta_{\chi_{2}}).
\]
In this case there is a Hilbert extension $\al F.$ of $\{\al A.,\Gamma\}$.

Conversely, if there is a Hilbert extension $\al F.$ of $\{\al
A.,\Gamma\}$, then to each $\alpha\in\Gamma$ there corresponds
a unitary $V_\alpha\in \bigcup_{\chi\in{\cal X}}{\cal
U}(\Pi_{\chi}{\cal F})$, such that $\alpha=\mbox{ad}\,V_\alpha\restriction
{\cal A}$ and
\[
 \epsilon(\alpha,\beta):=V_\alpha\,V_\beta\,V_\alpha^{-1}\,V_\beta^{-1}\,,
\]
is a system of representatives of the permutators 
$\widehat{\epsilon}(\alpha,\beta)$ satisfying the equations above. 
\end{pro}
\begin{rem}
\begin{itemize}
\item[(i)] In the present case the 2-cocycles $\lambda$ 
of the preceding section are $\T\EINS$-valued and
the relation (\ref{cocycle}) becomes the usual cocycle equation
\[
\lambda(\chi_{1},\chi_{2})\lambda(\chi_{1}\chi_{2},\chi_{3})=
\lambda(\chi_{2},\chi_{3})\lambda(\chi_{1},\chi_{2}\chi_{3}).
\]
\item[(ii)] In the particular case where $\al A.$ is the inductive
limit of a net of von Neumann algebras (which is a standard 
situation in algebraic quantum field theory, $\al A.$ being the
so--called quasilocal algebra) it can be shown that there is 
a lift $\gamma_\al X.$ of a given system of representatives 
$\beta_\al X.$, $\beta_\chi\in\chi$
(cf.~Definition~\ref{lift}), and by Corollary~\ref{1lift}
we have that
$\omega(\chi_1,\chi_2)=1$ is an admissible 2--cocycle of the 
system $\gamma_\al X.$. For a detailed construction of the 
lift see \cite{DHR69b}, \cite[Section~3.2]{bBaumgaertel95}.
\end{itemize}
\end{rem}

\section{A Hilbert space representation of $\{{\cal F},\alpha_{\cal G}\}$}
\label{Sec5}

Following Sutherland \cite{Sutherland80,SutherlandIn82}
one can introduce a faithful Hilbert
space representation of a Hilbert extension
$\{{\cal F},\alpha_{\cal G}\}$
of $\{{\cal A},\Theta({\cal X})\}.$

First let ${\cal H}$ be a Hilbert space and let $\pi$ be a
faithful representation of ${\cal A}$ on ${\cal H}$. Form the
Hilbert space
${\cal K}:=l^{2}({\cal X},{\cal H})$ 
by completion of
$C_{0}({\cal X}\rightarrow {\cal H})$
w.r.t. the norm
$\Vert f \Vert^{2}:=\sum_{\chi}\Vert f(\chi)\Vert^{2}_{\cal H}$.
Choose a system $\beta({\cal X})$ of representatives of
$\Theta({\cal X})$ and let $\omega$ be a corresponding generalized
2-cocycle such that
$U_{\chi_{1}}\cdot U_{\chi_{2}}=\omega(\chi_{1},\omega_{2})U_
{\chi_{1}\chi_{2}}.$
Now define a representation $\Phi$ of ${\cal F}_{0}\subset {\cal
F}$ on ${\cal K}$ by
\begin{eqnarray*}
(\Phi(A)f)(\chi)&:=&\pi(\beta_{\chi^{-1}}(A))f(\chi),\quad A\in {\cal A},\\
\Phi(U_{\chi_{0}})f)(\chi)
    &:=&\pi(\omega(\chi^{-1},\chi_{0}))f(
        \chi_{0}^{-1}\chi),\quad \chi_{0}\in {\cal X},\\
\Phi(AU_{\chi})
    &:=& \Phi(A)\Phi(U_{\chi}),\quad A\in{\cal A},\,\chi\in{\cal X}.
\end{eqnarray*}
Note that $\Phi(\EINS)=\EINS_{\cal K}$ and
$\Vert \Phi(A)\Vert_{\cal K}=\Vert A\Vert.$
One calculates easily
\begin{eqnarray*}
 \Phi(U_{\chi_{1}})\Phi(U_{\chi_{2}})
&=&\Phi(\omega(\chi_{1},\chi_{2}))\Phi(U_{\chi_{1}\chi_{2}}),\\
 \Phi(U_{\chi})\Phi(A)
&=&\Phi(\beta_{\chi}(A))\Phi(U_{\chi}),\\
 \Phi(A^{\ast})=\Phi(A)^{\ast},\quad
\Phi(U_{\chi}^{\ast})
&=&\Phi(U_{\chi})^{\ast}.
\end{eqnarray*}
Further $\Phi(\sum_{\chi}A_{\chi}U_{\chi})=0$ implies 
$\sum_{\chi}A_{\chi}U_{\chi}=0$,
i.e.~$\Phi$ is a *-isomorphism from ${\cal F}_{0}$ onto
$\Phi({\cal F}_{0})\subset {\cal L}({\cal K}).$
Recall that
\[
\Vert \Phi(F)\Vert_{\cal
K}=\mbox{sup}_{\Vert f\Vert\leq 1}\,\Vert\Phi(F)f\Vert_{\cal K}.
\]
We have
\begin{lem}
The relation
\begin{equation}\label{ineq}
\sup_{g\in{\cal G}}\Vert\Phi(\alpha_{g}F)\Vert_{\cal K}<\infty,\quad
F\in{\cal F}_{0},
\end{equation}
holds.
\end{lem}
\begin{beweis}
With $F=\sum_{\chi}A_{\chi}U_{\chi}$ we have
\begin{eqnarray*}
\Vert \Phi(F)f\Vert^{2}
&=&\sum_{y\in{\cal X}}\Vert \sum_{\chi}\pi(\alpha_{y^{-1}}(A_{\chi})
   \omega(y^{-1},\chi))f(y^{-1}\chi)\Vert_{\cal H}^{2}\\
&\leq&
  \sum_{y\in{\cal X}}(\sum_{\chi}\Vert \pi(\alpha_{y^{-1}}(A_{\chi})
  \omega(y^{-1},\chi))f(y^{-1}\chi)\Vert)^{2}\\
&\leq&
  \sum_{y\in{\cal X}}(\sum_{\chi}\Vert A_{\chi}\Vert \cdot
  \Vert f(y^{-1}\chi)\Vert)^{2}\leq
\sum_{y\in{\cal X}}(\sum_{\chi}\Vert A_{\chi}\Vert^{2})
(\sum_{\chi}\Vert f(y^{-1}\chi)\Vert^{2})\\
&=&
  (\sum_{\chi}\Vert A_{\chi}\Vert^{2})\sum_{\chi}
  \sum_{y\in{\cal X}}\Vert f(y^{-1}\chi)\Vert^{2} =
N(F)\Vert f\Vert^{2}\sum_{\chi}\Vert A_{\chi}\Vert^{2},
\end{eqnarray*}
where $N(F)$ denotes the number of terms of $F$. Hence we obtain
\[
\Vert \Phi(F)\Vert_{\cal K}\leq
N(F)^{1/2}(\sum_{\chi}\Vert A_{\chi}\Vert^{2})^{1/2}=:C_{F}.
\]
and this implies
\[
\Vert\Phi(\alpha_{g}F)\Vert_{\cal K}\leq C_{F},\quad g\in
{\cal G},
\]
because the number of terms of $\alpha_{g}F$ equals that of $F$
and $\Vert\chi(g)A_{\chi}\Vert = \Vert A_{\chi}\Vert.$
This implies the inequality (\ref{ineq}).
\end{beweis}

This result means that
\[
\Vert\Phi(F)\Vert_{sup}:=\sup_{g\in{\cal G}}
\Vert \Phi(\alpha_{g}F)\Vert_{\cal K}
\]
is a C*-norm on ${\cal F}_{0}.$

\begin{teo}\label{Teo3}
The relation
\[
\Vert\Phi(F)\Vert_{sup}=\Vert F\Vert_{\ast},\quad
F\in{\cal F}_{0}\,,
\]
holds, and in particular
$\Vert\Phi(F)\Vert_{\cal K}\leq \Vert F\Vert_{\ast}$, $F\in{\cal F}_{0}$.
\end{teo}
\begin{beweis}
The norm ${\cal F}_{0}\ni F\rightarrow \Vert
\Phi(F)\Vert_{sup}$ has the properties
$\Vert\Phi(A)\Vert_{sup}=\Vert A\Vert$ for all $A\in
{\cal A}$ and
$\Vert\Phi(\alpha_{g}F)\Vert_{sup}=\Vert\Phi(F)\Vert_
{sup}$ for all $g\in {\cal G}.$ However, according to
Doplicher/Roberts \cite[p.~105]{Doplicher87} there is at most one C*-norm on
${\cal F}_{0}$ with the mentioned properties.
\end{beweis}

\begin{rem}
If there is a faithful state $\phi_{0}$ of ${\cal A}$,
then Theorem~\ref{Teo3} can be improved. In this case
\[
\Vert \Phi(F)\Vert_{\cal K}=\Vert F \Vert_{\ast},\quad F\in
{\cal F}_{0},
\]
holds. This is implied by the fact that in this case
Sutherland's representation $\Phi$ of ${\cal F}_{0}$ on ${\cal
K}$ is unitarily equivalent to the so--called regular
representation of $\{{\cal F},\alpha_{\cal G}\}$ (restricted to 
${\cal F}_{0}$) given by the (faithful) GNS-representation
$\pi$ of $\{{\cal F},\alpha_{\cal G}\}$
on the GNS-Hilbert space ${\cal H}_{\pi}$ w.r.t. the ${\cal
G}$-invariant state
$\phi(F):=\phi_{0}(\Pi_{\iota}F),\,F\in{\cal F},$
such that
$\Vert\Phi(F)\Vert_{\cal K}=\Vert\pi(F)\Vert_{{\cal H}_{\pi}}$
for all $F\in {\cal F}_{0}$, but
$\Vert\pi(F)\Vert_{{\cal H}_{\pi}}=\Vert F\Vert_{\ast}$
for all $F\in {\cal F}$
(see, for example, \cite[p.~108 ff.]{bBaumgaertel95}).
\end{rem}

\paragraph{Acknowledgements}
One of us (H.B.) wants to thank Alan Carey for his kind invitation
to the University of Adelaide as well as for valuable suggestions
of an earlier version of the manuscript. The other author (F.Ll.)
expresses his gratitude to Sergio Doplicher for his hospitality
at the `Dipartamento di Matematica dell' Universit\`a di Roma `La
Sapienza'' in october '99. The visit was supported by a 
EU TMR network ``Implementation of concept and methods from
Non--Commutative Geometry to Operator Algebras and its applications'', 
contract no.~ERB FMRX-CT 96-0073.

\providecommand{\bysame}{\leavevmode\hbox to3em{\hrulefill}\thinspace}

\end{document}